\documentclass{amsart}

\usepackage{graphicx}
\usepackage{morefloats,bm,caption,amsbsy,enumerate,amsmath,amsthm,amssymb,mathtools,amsfonts,multirow,verbatim,tikz,tikz-cd,thmtools,thm-restate}
\usepackage{enumitem}
\usepackage[alphabetic]{amsrefs}
\usepackage{chngcntr}

\usepackage[hidelinks,colorlinks=true,linkcolor=blue, citecolor=black,linktocpage=true]{hyperref}
\usetikzlibrary{matrix,arrows,decorations.pathmorphing}
\usepackage[top=1.3in, bottom=1.1in, left=1.3in, right=1.3in,marginpar=1in]{geometry}
\usepackage{url}
\counterwithin{figure}{section}

\newtheorem{thm}{Theorem}[section]
\newtheorem{prop}[thm]{Proposition}

\newtheorem{lem}[thm]{Lemma}

\theoremstyle{definition}
\newtheorem{define}[thm]{Definition}

\theoremstyle{remark}

\newcommand{\ve}[1]{\boldsymbol{\mathbf{#1}}}
\newcommand{\R}{\mathbb{R}}

\newcommand{\Z}{\mathbb{Z}}

\renewcommand{\d}{\partial}
\renewcommand{\subset}{\subseteq}

\newcommand{\iso}{\cong}

\DeclareMathOperator{\Diff}{{Diff}}

\DeclareMathOperator{\ev}{{ev}}

\DeclareMathOperator{\id}{{id}}

\DeclareMathOperator{\Int}{{int}}

\DeclareMathOperator{\MCG}{{MCG}}

\DeclareMathOperator{\Spin}{{Spin}}

\DeclareMathOperator{\Sym}{{Sym}}

\DeclareMathOperator{\Tors}{{Tors}}

\DeclareMathOperator{\Cob}{\mathsf{Cob}}

\newcommand{\bF}{\mathbb{F}}

\newcommand{\bK}{\mathbb{K}}

\newcommand{\bT}{\mathbb{T}}

\newcommand{\cB}{\mathcal{B}}

\newcommand{\cH}{\mathcal{H}}

\newcommand{\cM}{\mathcal{M}}

\newcommand{\cW}{\mathcal{W}}

\newcommand{\frs}{\mathfrak{s}}

\newcommand{\CF}{\mathit{CF}}
\newcommand{\HF}{\mathit{HF}}

\newcommand{\HFh}{\widehat{\mathit{HF}}}
\newcommand{\CFh}{\widehat{\mathit{CF}}}

\newcommand{\SFH}{\mathit{SFH}}

\newcommand{\xs}{\ve{x}}
\newcommand{\ys}{\ve{y}}

\newcommand{\ps}{\ve{p}}
\newcommand{\as}{\ve{\alpha}}
\newcommand{\bs}{\ve{\beta}}

\renewcommand{\a}{\alpha}
\renewcommand{\b}{\beta}
\newcommand{\g}{\gamma}

\newcommand{\G}{\Gamma}

\usepackage{leftidx}

\renewcommand{\hat}{\widehat}

\newcommand{\Sp}{\mathit{Sp}}

\newcommand{\St}{S}
\newcommand{\M}{M}

\title{A graph TQFT for hat Heegaard Floer homology}
\author{Ian Zemke}
\address{Department of Mathematics\\Princeton University\\  Princeton, NJ 08544, USA}
\email{izemke@math.princeton.edu}
\thanks{This research was supported by NSF grant DMS-1703685.}
\begin{document}
\maketitle
\begin{abstract}
We construct maps on hat Heegaard Floer homology for cobordisms decorated with graphs. The graph TQFT allows for cobordisms with disconnected ends. Our construction uses Juh\'{a}sz's sutured Floer TQFT. We compute the maps for several elementary graph cobordisms. As an application, we compute the action of the fundamental group on hat Heegaard Floer homology.
\end{abstract}

\section{Introduction}

Ozsv\'{a}th and Szab\'{o} constructed a powerful set of invariants for closed 3-manifolds, and cobordisms between them \cite{OSDisks} \cite{OSTriangles}. To a closed, oriented 3-manifold $Y$, they constructed a finitely generated abelian group
\[
\hat{\HF}(Y),
\]
as well as $\Z[U]$-modules $\HF^-(Y)$, $\HF^+(Y)$ and $\HF^\infty(Y)$. We focus mostly on $\hat{\HF}$ in our present paper. Also, we work over the field $\bF:=\Z/2\Z$.

To a compact, connected, and oriented cobordism $W$ between two connected 3-manifolds $Y_0$ to $Y_1$, they constructed a linear map
\[
\hat{F}_{W}\colon \hat{\HF}(Y_1)\to \hat{\HF}(Y_2).
\]
If $W=W_2\cup_Y W_1$, where $Y$ is a closed, connected 3-manifold, then 
\[
\hat{F}_{W}=\hat{F}_{W_2}\circ \hat{F}_{W_1}.
\]

An important subtlety is that the construction of $\hat{\HF}(Y)$ requires a choice of basepoint in $Y$. Similarly, the construction of $\hat{F}_{W}$ implicitly relies on choosing an arc in $W$, connecting the two basepoints in $\d W$.  To make the dependence explicit, we will write $\hat{\HF}(Y,p)$ and $\hat{F}_{W,\gamma}$, for the groups and maps defined with an auxiliary choice of basepoint or arc.

\subsection{Maps for graph cobordisms}

The main construction of this paper is an extension of Ozsv\'{a}th and Szab\'{o}'s cobordism maps to the following category:

\begin{define}The \emph{graph cobordism category} $\Cob^\G_{3+1}$ has the following objects and morphisms:
 \begin{itemize}
\item The objects are pairs $(Y,\ve{p})$, where $Y$ is a closed and oriented 3-manifold (possibly disconnected or empty), and $\ve{p}$ is a finite collection of basepoints in $Y$, such that each component of $Y$ has at least one basepoint.
\item A morphism from $(Y_0,\ve{p}_0)$ to $(Y_1,\ve{p}_1)$ is a pair $(W,\Gamma)$ such that
\begin{enumerate}
\item $W$ is an oriented, compact cobordism from $Y_0$ to $Y_1$, and
\item $\Gamma\subset W$ is an embedded graph, such that $\Gamma\cap Y_i=\ve{p}_i$, $\Gamma$ has no valence 0 vertices,  and $\ve{p}_i\subset \Gamma$ are all valence 1.
\end{enumerate}
\end{itemize}
\end{define}

Generalizing their construction of Heegaard Floer homology for singly pointed 3-manifolds \cite{OSDisks}, Ozsv\'{a}th and Szab\'{o} also 
defined a group $\hat{\HF}(Y,\ve{p})$, whenever $(Y,\ve{p})$ is a closed, oriented 3-manifold with a collection of basepoints \cite{OSLinks}. The construction extends via a tensor product to disconnected 3-manifolds, as long as each component of $Y$ contains at least one basepoint.

In this paper, we construct cobordism maps for graph cobordisms, and prove the following:

\begin{thm}\label{thm:graph-cobordism-maps} If $(W,\Gamma)$ is a graph cobordism from $(Y_0,\ve{p}_0)$ to $(Y_1,\ve{p}_1)$, then the construction of this paper gives a well-defined map
\[
\hat{F}_{W,\Gamma}\colon \hat{\HF}(Y_0,\ve{p}_0)\to \hat{\HF}(Y_1,\ve{p}_1),
\]
satisfying the following:
\begin{enumerate}
\item  $\hat{F}_{[0,1]\times Y, [0,1]\times \ve{p}}=\id_{\hat{\HF}(Y,\ve{p})}$.
\item If $(W,\Gamma)=(W_2,\Gamma_2)\cup (W_1,\Gamma_1)$, then
\[
\hat{F}_{W,\Gamma}=\hat{F}_{W_2,\Gamma_2}\circ \hat{F}_{W_1,\Gamma_1}.
\]
\item If $(W,\Gamma)\colon (Y_0,p_1)\to (Y_1,p_1)$ is a cobordism such that $\Gamma$ is a path connecting $p_0$ and $p_1$, then $\hat{F}_{W,\Gamma}$ coincides with the map of Ozsv\'{a}th and Szab\'{o} \cite{OSTriangles}.
\end{enumerate}
\end{thm}

Theorem~\ref{thm:graph-cobordism-maps} implies that Heegaard Floer homology extends to a functor from $\Cob_{3+1}^\G$ to the category of vector spaces over $\bF$. Our construction of $\hat{F}_{W,\Gamma}$ uses Juh\'{a}sz's sutured Floer homology TQFT \cite{JDisks} \cite{JCob}.

\subsection{Elementary graph cobordisms}

In this paper,  we  compute the maps for the following elementary graph cobordisms. See Sections ~\ref{sec:elementary-computations} and~\ref{sec:broken-path}.
\begin{enumerate}[label=($\Gamma$-\arabic*), ref=$\Gamma$-\arabic*]
\item\label{item:free-stab}  \emph{Free-stabilization cobordisms}:  The graph $\Gamma$ consists of $n\ge 1$ arcs connecting $\{0\}\times Y$ to $\{1\}\times Y$, as well as one half-arc of the form $[0,\tfrac{1}{2}]\times \{p\}$ or $[\tfrac{1}{2},1]\times\{p\}$, for some $p\in Y$.
\item\label{item:split/merge} \emph{Merging and splitting cobordisms}: The graph consists of $n\ge 0$ strands from $\{0\}\times Y$ to $\{1\}\times Y$, as well as one wye-shaped component which merges or splits two basepoints.
\item\label{item:loop-splicing} \emph{Spliced loop cobordisms}: The graph consists of $n\ge 1$ strands going from $\{0\}\times Y$ to $\{1\}\times Y$, as well as one loop $\gamma$ in $\{\tfrac{1}{2}\}\times Y$, which intersects one of the vertical strands at a single point.
\item\label{item:broken-path} \emph{Broken path cobordisms}: The cobordism consists of $n\ge 0$ arcs going from $\{0\}\times Y$ to $\{1\}\times Y$, as well as one broken arc $([0,\tfrac{1}{3}]\cup [\tfrac{2}{3},1])\times \{p\}$.
\end{enumerate}

\begin{figure}[ht!]
	\centering
%% Accompanies image file 'fig1.pdf' (pdf, eps, ps)
\begingroup%
  \makeatletter%
  \providecommand\color[2][]{%
    \errmessage{(Inkscape) Color is used for the text in Inkscape, but the package 'color.sty' is not loaded}%
    \renewcommand\color[2][]{}%
  }%
  \providecommand\transparent[1]{%
    \errmessage{(Inkscape) Transparency is used (non-zero) for the text in Inkscape, but the package 'transparent.sty' is not loaded}%
    \renewcommand\transparent[1]{}%
  }%
  \providecommand\rotatebox[2]{#2}%
  \newcommand*\fsize{\dimexpr\f@size pt\relax}%
  \newcommand*\lineheight[1]{\fontsize{\fsize}{#1\fsize}\selectfont}%
  \ifx\svgwidth\undefined%
    \setlength{\unitlength}{171.7052241bp}%
    \ifx\svgscale\undefined%
      \relax%
    \else%
      \setlength{\unitlength}{\unitlength * \real{\svgscale}}%
    \fi%
  \else%
    \setlength{\unitlength}{\svgwidth}%
  \fi%
  \global\let\svgwidth\undefined%
  \global\let\svgscale\undefined%
  \makeatother%
  \begin{picture}(1,0.97900609)%
    \lineheight{1}%
    \setlength\tabcolsep{0pt}%
    \put(0,0){\includegraphics[width=\unitlength,page=1]{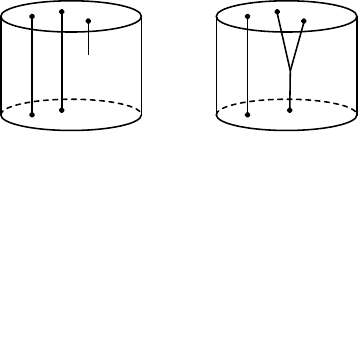}}%
    \put(0.19920432,0.54657897){\color[rgb]{0,0,0}\makebox(0,0)[t]{\lineheight{1.25}\smash{\begin{tabular}[t]{c}\eqref{item:free-stab}\end{tabular}}}}%
    \put(0,0){\includegraphics[width=\unitlength,page=2]{fig1.pdf}}%
    \put(0.80172033,0.5468292){\color[rgb]{0,0,0}\makebox(0,0)[t]{\lineheight{1.25}\smash{\begin{tabular}[t]{c}\eqref{item:split/merge}\end{tabular}}}}%
    \put(0.19920432,0.00962314){\color[rgb]{0,0,0}\makebox(0,0)[t]{\lineheight{1.25}\smash{\begin{tabular}[t]{c}\eqref{item:loop-splicing}\end{tabular}}}}%
    \put(0,0){\includegraphics[width=\unitlength,page=3]{fig1.pdf}}%
    \put(0.80172033,0.00937291){\color[rgb]{0,0,0}\makebox(0,0)[t]{\lineheight{1.25}\smash{\begin{tabular}[t]{c}\eqref{item:broken-path}\end{tabular}}}}%
    \put(0.26130433,0.28259954){\color[rgb]{0,0,0}\makebox(0,0)[lt]{\lineheight{1.25}\smash{\begin{tabular}[t]{l}$\g$\end{tabular}}}}%
  \end{picture}%
\endgroup%

	\caption{The four elementary graph cobordisms in $[0,1]\times Y$.}\label{fig:1}
\end{figure}

The  elementary graph cobordisms \eqref{item:free-stab}--\eqref{item:broken-path} are depicted in Figure~\ref{fig:1}.

\subsection{The action of the fundamental group}

Since a basepoint is implicitly used in the construction of the Heegaard Floer groups, the naturality theorem of \cite{JTNaturality}
implies only that elements of the \emph{based} mapping class group $\MCG(Y,p)$ induce well defined endomorphisms of Heegaard Floer homology.

There is a fibration
\[
\Diff(Y,p)\to \Diff(Y)\xrightarrow{\ev_p} Y
\]
where $\ev_p$ denotes evaluation at $p$. The long exact sequence of homotopy groups for this fibration gives a homomorphism
\[
\pi_1(Y,p)\to \MCG(Y,p).
\]
By exactness, the image of $\pi_1(Y,p)$ in $\MCG(Y,p)$ is the kernel of the forgetful map $\MCG(Y,p)\to \MCG(Y)$.

Suppose $p\in \ve{p}$ and $\g\in \pi_1(Y,p)$. We write
\[
\gamma_*\colon \hat{\HF}(Y,\ve{p})\to \hat{\HF}(Y,\ve{p})
\]
for the induced endomorphism.

Using  the graph TQFT, we prove the following:
\begin{thm}
Suppose $(Y,\ve{p})$ is a multi-pointed 3-manifold and $p\in \ve{p}$. If $\gamma\in \pi_1(Y,p)$, then the induced endomorphism $\gamma_*$ of $\hat{\HF}(Y,\ve{p})$ satisfies
 \[
 \gamma_*=\id+\Phi_p \circ A_{[\gamma]},
 \]
 where $A_{[\gamma]}$ denotes the standard action of  $[\gamma]\in H_1(Y;\Z)/\Tors$, and $\Phi_p\colon \hat{\HF}(Y,\ve{p})\to \hat{\HF}(Y,\ve{p})$ is the broken path graph cobordism labeled \eqref{item:broken-path} in Figure~\ref{fig:1}.
\end{thm}

In Section~\ref{sec:broken-path}, we identify the broken path graph cobordism $\Phi_p\colon \HFh(Y,\ve{p})\to \HFh(Y,\ve{p})$ with the basepoint action for the point $p$, which counts holomorphic disks on a Heegaard diagram with multiplicity 1 at $p$. See Proposition~\ref{prop:d1-comp}.

\subsection*{Acknowledgments}

I would like to thank  Jianfeng Lin, Yajing Liu, Andr\'{a}s Juh\'{a}sz, Ko Honda, Ciprian Manolescu, Marco Marengon  and Matthew Stoffregen for helpful conversations.

\section{Background}

\subsection{Heegaard Floer homology}  Suppose $(Y,\ve{p})$ is a multi-pointed 3-manifold, $\frs\in \Spin^c(Y)$, and $\cH=(\Sigma,\as,\bs,\ps)$ is a Heegaard diagram for $(Y,\ve{p})$. Ozsv\'{a}th and Szab\'{o} \cite{OSLinks} construct chain complexes $\CFh(\cH,\frs)$, $\CF^-(\cH,\frs)$, $\CF^+(\cH,\frs)$ and $\CF^\infty(\cH,\frs)$, as follows. We focus on the case that $Y$ is connected. If $Y$ is disconnected, then $\CFh(\cH,\frs)$ is defined by tensoring over $\bF$ the complexes for each connected component.

 The chain complex $\CFh(\cH,\frs)$ is generated by intersection points $\xs$ of the two half dimensional tori
\[
\bT_{\a}=\a_1\times \cdots \times \a_{n+g(\Sigma)-1}\quad \text{and} \quad \bT_{\b}=\b_1\times \cdots \times \b_{n+g(\Sigma)-1},
\]
in $\Sym^{g(\Sigma)+n-1}(\Sigma)$ (where $n=|\ps|$), which satisfy $\frs_{\ps}(\xs)=\frs$. The differential is given by the formula
\[
\d \xs=\sum_{\ys\in \bT_{\a}\cap \bT_{\b}} \sum_{\substack{\phi\in \pi_2(\xs,\ys)\\ \mu(\phi)=1\\ n_{\ve{p}}(\phi)=0}} \# (\cM(\phi)/\R) \cdot \ys,
\] 
where $\cM(\phi)$ denotes the moduli space of holomorphic disks in $\Sym^{g(\Sigma)+n-1}(\Sigma)$ representing a given homotopy class $\phi\in \pi_2(\xs,\ys)$.

We define
\begin{equation}
\CFh(\cH)=\bigoplus_{\frs\in \Spin^c(Y)} \CFh(\cH,\frs). \label{eq:CFh-direct-sum}
\end{equation}

Although we mostly focus on $\CFh$ in this paper, in Section~\ref{sec:broken-path}, we consider $\CF^-$, which we review presently. Write $\ve{p}=\{p_1,\dots, p_n\}$, and $R_n:=\bF[U_1,\dots, U_n]$. The module $\CF^-(\cH,\frs)$ is the free $R_n$-module with generators $\xs\in \bT_{\a}\cap \bT_{\b}$ with $\frs_{\ve{p}}(\xs)=\frs$. The differential on $\CF^-(\cH,\frs)$ is
\[
\d \xs= \sum_{\ys\in \bT_{\a}\cap \bT_{\b}} \sum_{\substack{\phi\in \pi_2(\xs,\ys)\\ \mu(\phi)=1}} \# (\cM(\phi)/\R)U_1^{n_{p_1}(\phi)}\cdots U_n^{n_{p_n}(\phi)} \cdot \ys.
\]

Unlike for $\CFh$, it is usually not possible to define a total complex $\CF^-(\cH)$ as a direct sum over all $\Spin^c$ structures, analogous to equation~\eqref{eq:CFh-direct-sum}, since $\CF^-$ requires a stronger version of admissibility than $\CFh$, which cannot normally be simultaneously realized for all $\Spin^c$ structures on a single diagram \cite{OSDisks}*{Section~4.2.2}. Hence, on $\CF^-$, we usually work with one $\Spin^c$ structure at a time.

\subsection{Sutured Floer homology}

 \emph{Sutured manifolds} were defined by Gabai \cite{GabaiFoliationsI} to study foliations on 3-manifolds.  Juh\'{a}sz \cite{JCob} defined an extension of Heegaard Floer homology for sutured manifolds, called \emph{sutured Floer homology}. Juh\'{a}sz \cite{JDisks} also described a $(3+1)$-dimensional TQFT for sutured Floer homology. In this section, we recall some background about sutured manifolds and the sutured Floer homology TQFT.
 
 The following is a slight restriction of Gabai's original definition, but is sufficient for our purposes:

\begin{define}\label{def:sutured-manifold}
 A \emph{sutured manifold} $(M,\g)$ is a compact, oriented 3-manifold $M$ with boundary, together with a set of pairwise disjoint, simple closed curves $\gamma\subset \d M$ (the sutures) which are oriented. The surface $\d M\setminus \g$ is partitioned into two sets of components, $R_+(\g)$ and $R_-(\g)$, which are oriented so that the normal of $R_+(\g)$ points out of $M$, while the normal of $R_-(\g)$ points into $M$. Finally, we require $\g$ to be consistently oriented with respect to the boundary orientation of $R_+(\g)$ and $R_-(\g)$.
\end{define}

The main difference between Definition~\ref{def:sutured-manifold} and Gabai's original definition is that we exclude toroidal sutures. We say that a sutured manifold $(M,\g)$ is \emph{balanced} if $\chi(R_+(\g))=\chi(R_-(\g))$.

To a balanced sutured manifold $(M,\g)$ with no closed components, Juh\'{a}sz constructs an $\bF$-vector space
\[
\SFH(M,\g).
\]

If $Y$ is a closed, oriented 3-manifold, and $\ve{p}$ is a collection of basepoints, then we write $Y(\ve{p})$ for the sutured manifold $(M,\g)$ where
\[
M:=Y\setminus \Int N(\ve{p})
\]
and $\g$ consists of one simple closed curve per component of $\d M$. We note that
\[
\SFH(Y(\ve{p}))=\hat{\HF}(Y,\ve{p}),
\]
since a Heegaard diagram for $(Y,\ve{p})$ may be obtained from a diagram for $Y(\ve{p})$ by collapsing each suture to a basepoint.

Juh\'{a}sz also defines cobordism maps for sutured Floer homology \cite{JCob}.
He uses the following notion of cobordism between sutured manifolds:
\begin{define}\label{def:sutured-cobordism} A \emph{cobordism of sutured
manifolds} 
	\[
	\cW=(W,Z,[\xi])\colon (M_0,\g_0)\to (M_1,\g_1)	
	\]
	is a triple such that
	\begin{enumerate}
		\item $W$ is a compact, oriented 4-manifold with boundary,
		\item $Z$ is a compact, codimension-0 submanifold
		 of $\d W$, and $\d W\setminus \Int (Z)=-M_0\sqcup M_1$,
		\item $[\xi]$ is an equivalence class of positive contact structures
		 on $Z$, such that $\d Z$ is a convex surface with dividing
		  set $\gamma_i$ on $\d M_i$, for $i\in \{0,1\}$.
	\end{enumerate}
	\end{define}
	The notion of equivalence between contact structures used in
	Definition~\ref{def:sutured-cobordism} can be found in 
	\cite{JCob}*{Definition~2.3}.

If $\cW\colon (M_0,\g_0)\to (M_1,\g_1)$ is a cobordism between balanced sutured manifolds, Juh\'{a}sz \cite{JCob} constructs  a well-defined map
\[
F_{\cW}\colon \SFH(M_0,\g_0)\to \SFH(M_1,\g_1),
\]
which is functorial in the following sense. If $\xi$ is a $[0,1]$-invariant contact structure on $[0,1]\times \d M$, such that $\{0,1\}\times \d M$ is convex, with dividing set $\g$, then $\cW=([0,1]\times M, [0,1]\times \d M,[\xi])$ is a sutured manifold cobordism from $(M,\g)$ to itself. The induced cobordism map
\[
F_{\cW}\colon \SFH(M,\g)\to \SFH(M,\g)
\]
is the identity. Furthermore, if $\cW$ decomposes as the composition of two sutured manifold cobordisms $\cW_2\circ \cW_1$, then
\[
F_{\cW}=F_{\cW_2}\circ F_{\cW_1}.
\]
See \cite{JCob}*{Theorem~11.3}.

We outline the construction of the sutured cobordism maps in Section~\ref{sec:outline-sutured-construction}, after we outline one of its constituents, the \emph{contact gluing map}.

\subsection{The contact gluing map}
\label{sec:gluing-map}
We now recall the Honda--Kazez--Mati\'{c} \emph{contact gluing map} for sutured Floer homology, as well as the contact-handle formulation given by Juh\'{a}sz and the author \cite{JuhaszZemkeContactHandles}.

\begin{define}
 Suppose that $(M,\g)$ and $(M',\g')$ are sutured manifolds. We say that $(M,\g)$ is a \emph{sutured submanifold} of $(M',\g')$ if $M\subset \Int(M')$.
\end{define}

If $(M,\g)$ is a sutured submanifold of $(M',\g')$, and $\xi$ is a positive contact structure on $Z:=M'\setminus \Int (M)$ which induces the dividing set $\g\cup \g'$, then Honda, Kazez and Mati\'{c} \cite{HKMTQFT} define a \emph{contact gluing map}
\[
\Phi_{Z,\xi}\colon \SFH(-M,\g)\to \SFH(-M',\g').
\]

In \cite{JuhaszZemkeContactHandles}, Juh\'{a}sz and the author gave a reformulation of the contact gluing map in terms of \emph{contact handles}, which facilitates computations. Contact handles were defined by Giroux \cite{Giroux}. See Ozbagci \cite{Ozbagci} for an exposition. We take the following as the definition of a contact handle:

\begin{define} 
Suppose $(M,\g)$ is a sutured submanifold of $(M',\g')$, and $\xi$ is a positive contact structure on $Z:=M'\setminus \Int(M)$, with dividing set $\g\cup \g'$. We say that $(Z,\xi)$ is a \emph{contact handle of index $k$} if there is a contact vector field $\nu$ on $Z$ that points into $Z$ on $\d M$, and out of $Z$ on $\d M'$, as well as a decomposition $Z=Z_0\cup h$ such that
\begin{enumerate}
\item $Z_0$ is diffeomorphic to $[0,1]\times \d M$,
\item $\nu$ is non-vanishing on $Z_0$, points into $Z_0$ on $\{0\}\times \d M$ and out of $Z_0$ on $\{1\}\times \d M$, and each flowline of $\nu$ is an arc from $\{0\}\times \d M$ to $\{1\}\times \d M$,
\item $h$ is smooth 3-ball with corners, and $\xi$ is tight on $h$.
\end{enumerate}
We have the following additional requirements, depending on $k$:
\end{define}
\begin{itemize}
\item $(k=0)$: $h=D^3$ (with no corners) and $h\cap Z_0=\emptyset$. The dividing set on $\d h$ is a single circle.
\item $(k=1)$:  $h=[0,1]\times D^2$, and $h\cap Z_0=\{0,1\}\times D^2$. The dividing set on $\d h$ is a single closed curve, consisting of an arc on $\{0\}\times D^2$ and $\{1\}\times D^2$, and two longitudinal arcs on $[0,1]\times \d D^2$.
\item $(k=2)$: $h=[0,1]\times D^2$, and $h\cap Z_0=[0,1]\times \d D^2$. The dividing set is as in $k=1$ case.
\item $(k=3)$: $h= D^3$ (with no corners), and $h\cap Z_0=\d h$. The dividing set on $\d h$ is a single circle.
\end{itemize}

We now state the description from \cite{JuhaszZemkeContactHandles} of the contact gluing maps of Honda, Kazez and Mati\'{c} when $M'\setminus M$ is a contact handle.

If $Z$ is a contact 0-handle, we extend the Heegaard surface into $Z_0$ using the flow of $\nu$, and then add a disk to the Heegaard surface which lies in $h$ and intersects $\d D^3$ along the sutures. We add no new alpha or beta curves. The map on sutured Floer homology is the tautological one. 

If $Z$ is a contact 1-handle, we extend $\Sigma$ over $Z_0$ using the flow of $\nu$, and then attach a band to the boundary of the Heegaard surface, which is contained in $h$ and intersects the boundary along the dividing set. We add no alpha or beta curves. Similar to the contact 0-handle map, the map on sutured Floer homology is the tautological one.

\begin{figure}[ht!]
	\centering
	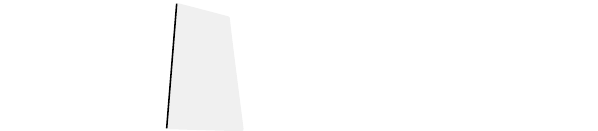
	\caption{A contact 2-handle on Heegaard diagrams.}\label{fig:9}
\end{figure}

If $Z$ is a contact 2-handle, we extend the Heegaard surface into $Z_0$ using the flow of $\nu$, as before. Now $\d h$ intersects $\d Z$ in an annulus. Let $c$ denote the core of the attaching annulus. The curve $c$ may be taken to intersect the dividing set in two points. Let $\lambda_+$ denote the subarc of $c$ which intersects $R_+$, and let $\lambda_-$ denote the subarc which intersects $R_-$. If $(\Sigma,\as,\bs)$ is a diagram for $(M,\g)$, we may obtain a diagram for $(M',\g')$ by adjoining a band to $\d \Sigma$, and adding a new alpha curve $\a$, and a new beta curve $\b$. The curves $\a$ and $\b$ have a single intersection point in the band region, as in Figure~\ref{fig:9}. Outside of the band region, $\b$ consists of $\lambda_+$, projected onto $\Sigma\setminus \bs$, and $\a$ consists of $\lambda_-$, projected on $\Sigma\setminus \as$. The map 
\[
\Phi_{Z,\xi}\colon \SFH(\Sigma,\bs,\as)\to \SFH(\Sigma\cup B, \bs\cup \{\b\}, \as\cup \{\a\})
\]
 is given by $\xs\mapsto \xs\times c$.

Finally, suppose $Z$ is a contact 3-handle, and let $S\subset \d M$ denote the 2-sphere which is filled in by $Z$. Let $S'$ denote a 2-sphere in $\Int (M)$ obtained by pushing $S$ into $\Int(M)$. The contact 3-handle map is defined as the composition of the 4-dimensional 3-handle map for the 2-sphere $S'$ (which leaves the disjoint union of $(M',\g')$ and $B^3$), followed by the canonical isomorphism
\[
\SFH(M',\g')\otimes \SFH(B^3)\iso \SFH(M',\g').
\]

\subsection{Sutured cobordism maps}
\label{sec:outline-sutured-construction}

We now outline the construction of the sutured cobordism maps. Suppose 
\[
\cW=(W,Z,[\xi])\colon (M_0,\g_0)\to (M_1,\g_1),
\]
is a cobordism of sutured manifolds, as in Definition~\ref{def:sutured-cobordism}.
 First, we remove some number of tight, contact 3-balls from $Z$, and add them to $(M_0,\g_0)$ or $(M_1,\g_1)$, so that each component of $W$ intersects a component of $Y_0$ and $Y_1$ non-trivially. This does not affect the sutured Floer homology of $(M_0,\g_0)$ or $(M_1,\g_1)$, as there is a canonical isomorphism
\[
\SFH(M_0\sqcup B^3, \g_0\cup \g)\iso \SFH(M_0,\g_0),
\]
where $\g\subset B^3$ denotes a single closed curve.

Juh\'{a}sz calls a sutured cobordism $\cW$ \emph{special} if $Z=[0,1]\times \d M_0$ and $\xi$ is $[0,1]$-invariant. The cobordism map for a special cobordism is constructed to be a composition of 1-handle, 2-handle and 3-handle maps, similar to the ones described by Ozsv\'{a}th and Szab\'{o} \cite{OSTriangles}.

If 
\[
\cW=(W,Z,[\xi])\colon (M_0,\g_0)\to (M_1,\g_1)
\]
is a general sutured manifold cobordism, one may obtain a special cobordism
\[
\cW^s=(W,[0,1]\times \d M_1, \xi_1)\colon (M_0\cup Z,\g_1)\to (M_1,\g_1),
\]
where $\xi_1$ is a $[0,1]$-invariant contact structure on $\d M_1\times [0,1]$.  The cobordism map $F_{\cW}$ is defined as the composition
\begin{equation}
F_{\cW}:=F_{\cW^s}\circ \Phi_{Z,\xi}.\label{eq:sutured-cobordism-definition}
\end{equation}

\section{Construction of the graph TQFT}

Suppose $(W,\Gamma)$ is a graph cobordism from $(Y_0,\ve{p}_0)$ to $(Y_1,\ve{p}_1)$. We define a sutured manifold cobordism
\[
\cW(W,\Gamma)=(W(\G),Z(\G),[\xi(\G)])\colon Y_0(\ve{p}_0)\to Y_1(\ve{p}_1),
\]
as follows. We define the 4-manifold $W(\G)$ to be $W\setminus \Int  N(\Gamma)$, and the 3-manifold $Z(\Gamma)$ to be $\d W(\G)\cap \d N(\Gamma)$. We give $\d Z(\G)$ the same sutures as $Y(\ve{p}_0)$ and $Y(\ve{p}_1)$, for which we write $\g_Z$. We take $\xi(\G)$ to be the unique tight contact structure with dividing set $\g_Z$, whose well definedness we sketch presently. The 3-manifold $Z(\G)$ is homeomorphic to a disjoint union of connected sums of $S^1\times S^2$, with some number of 3-balls removed. The sutures consist of a single closed curve on each copy of $S^2$ in $\d Z(\G)$.  It is well known that up to isotopy, relative to $\d Z(\G)$, there is a unique tight contact structure on $Z(\G)$ which has this dividing set on $\d Z(\G)$. The case when $Z(\G)=B^3$ follows from a result of Eliashberg \cite{Eliashberg-tight-S3}. The general case follows by decomposing $Z(\G)$ along a collection of convex 2-spheres, until one obtains a disjoint union of tight, punctured 3-spheres, using convex surface theory \cite{ColinGluing} \cite{HondaGluing}.

Without further ado, we define
\[
\hat{F}_{W,\G}:=F_{\cW(W,\G)}.
\]

\section{Elementary graph cobordisms in $[0,1]\times Y$}
\label{sec:elementary-computations}

In this section, we compute the maps induced by the elementary graph cobordisms shown in Figure~\ref{fig:1}, with the exception of the broken path cobordism, which is considered in Section~\ref{sec:broken-path}.

\subsection{Free-stabilization cobordisms}
\label{sec:free-stabilization}
In this section we compute the maps for the free-stabilization cobordisms, which are labeled \eqref{item:free-stab} in Figure~\ref{fig:1}.  Let 
\[
\cW_{p}^+:=([0,1]\times Y,\G^+_p)\colon (Y,\ps)\to (Y,\ps\cup \{p\})
\]
 denote the free-stabilization graph cobordism which adds the basepoint $p$, and let $\cW_p^-$ denote the free-stabilization graph cobordism which removes $p$.
 
 We define
\begin{equation}
\St_{p}^+:=\hat{F}_{\cW_p^+}\quad \text{and} \quad \St_p^-:=\hat{F}_{\cW_{p}^-}.
\end{equation}

\begin{figure}[ht!]
	\centering
%% Accompanies image file 'fig2.pdf' (pdf, eps, ps)
\begingroup%
  \makeatletter%
  \providecommand\color[2][]{%
    \errmessage{(Inkscape) Color is used for the text in Inkscape, but the package 'color.sty' is not loaded}%
    \renewcommand\color[2][]{}%
  }%
  \providecommand\transparent[1]{%
    \errmessage{(Inkscape) Transparency is used (non-zero) for the text in Inkscape, but the package 'transparent.sty' is not loaded}%
    \renewcommand\transparent[1]{}%
  }%
  \providecommand\rotatebox[2]{#2}%
  \newcommand*\fsize{\dimexpr\f@size pt\relax}%
  \newcommand*\lineheight[1]{\fontsize{\fsize}{#1\fsize}\selectfont}%
  \ifx\svgwidth\undefined%
    \setlength{\unitlength}{101.93527793bp}%
    \ifx\svgscale\undefined%
      \relax%
    \else%
      \setlength{\unitlength}{\unitlength * \real{\svgscale}}%
    \fi%
  \else%
    \setlength{\unitlength}{\svgwidth}%
  \fi%
  \global\let\svgwidth\undefined%
  \global\let\svgscale\undefined%
  \makeatother%
  \begin{picture}(1,0.67569496)%
    \lineheight{1}%
    \setlength\tabcolsep{0pt}%
    \put(0,0){\includegraphics[width=\unitlength,page=1]{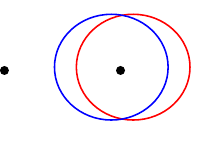}}%
    \put(0.82763266,0.55929479){\color[rgb]{1,0,0}\makebox(0,0)[lt]{\lineheight{1.25}\smash{\begin{tabular}[t]{l}$\alpha$\end{tabular}}}}%
    \put(0.32089156,0.55831779){\color[rgb]{0,0,1}\makebox(0,0)[rt]{\lineheight{1.25}\smash{\begin{tabular}[t]{r}$\beta$\end{tabular}}}}%
    \put(0.59673053,0.00691692){\color[rgb]{0,0,0}\makebox(0,0)[t]{\lineheight{1.25}\smash{\begin{tabular}[t]{c}$\theta^+$\end{tabular}}}}%
    \put(0.58774458,0.64549811){\color[rgb]{0,0,0}\makebox(0,0)[t]{\lineheight{1.25}\smash{\begin{tabular}[t]{c}$\theta^-$\end{tabular}}}}%
    \put(0.60837962,0.38908337){\color[rgb]{0,0,0}\makebox(0,0)[lt]{\lineheight{1.25}\smash{\begin{tabular}[t]{l}$p$\end{tabular}}}}%
    \put(0.06396307,0.38908337){\color[rgb]{0,0,0}\makebox(0,0)[lt]{\lineheight{1.25}\smash{\begin{tabular}[t]{l}$p_0$\end{tabular}}}}%
  \end{picture}%
\endgroup%

	\caption{The diagram $\cH_{(p)}$, obtained by adding a basepoint $p$ to a diagram $\cH$.}\label{fig:2}
\end{figure}

If $\cH$ is a Heegaard diagram for $(Y,\ps)$, we may form a Heegaard diagram $\cH_{(p)}$ for $(Y,\ve{p}\cup \{p\})$ by adding the basepoint $p$, encircled by a new pair of alpha and beta curves, $\a$ and $\b$, as in Figure~\ref{fig:2}. By picking $\cH$ appropriately, we may assume that $\a$ and $\b$ are immediately adjacent to another basepoint $p_0\in \ve{p}$. The choice of basepoints makes it easy to verify that
\begin{equation}
\HFh(\cH_{(p)})\iso \HFh(\cH)\otimes V,\label{eq:isomorphism-extra-basepoint}
\end{equation}
where $V$ is the 2-dimensional vector space $\bF_{1/2}\oplus \bF_{-1/2}$. We write $\theta^+$ for the top degree generator of $V$, and $\theta^-$ for the bottom degree generator. 

\begin{lem}\label{lem:compute-free-stabilization}With respect to the isomorphism from equation~\eqref{eq:isomorphism-extra-basepoint}, the maps $S_p^+$ and $S_p^-$ satisfy
\[
\St_p^+(\xs)=\xs\times \theta^+\quad \text{and} \quad \St_p^-(\xs\times \theta)=\begin{cases} \xs& \text{if } \theta=\theta^-\\
0& \text{if } \theta=\theta^+.
\end{cases}
\]
\end{lem}
\begin{proof} We first consider $\St^+_p$. We may perform an index 0/1 handle cancellation to decompose the graph cobordism $\cW_p^+$ as follows:
\begin{enumerate}
\item A 0-handle $B^4$, containing an arc $a$, which intersects $\d B^4$ in a single point.
\item A 1-handle cobordism which merges $(S^3,p)$ with $(Y,\ve{p})$, away from $\ve{p}$.
\end{enumerate}
We may similarly decompose $\cW_p^-$ into a 3-handle cobordism followed by a 4-handle cobordism.

The graph cobordism map for $(B^4,a)\colon \emptyset \to (S^3,p)$ is easily seen to send the generator of $\HFh(\emptyset)\iso \bF$ to the generator of $\HFh(S^3)\iso \bF$, and similarly for the 4-handle cobordism in the opposite direction.   The main claim now follows for $S_{p}^+$, since the stated formula coincides with the definition of the 1-handle map \cite{JCob}*{Section~7}. The proof of $S_p^-$ is similar.
\end{proof}

\subsection{Merge and split cobordisms}

We now compute the merge and split cobordism maps, which are labeled~\eqref{item:split/merge} in Figure~\ref{fig:1}. Suppose that $p_1$ and $p_2$ are two points in $Y$, $\lambda$ is a path connecting $p_1$ to $p_2$, and $p_0$ is a point along $\lambda$. Suppose that $\ve{p}$ is a (possibly empty) collection of basepoints in $Y\setminus\{p_0,p_1,p_2\}$. Write 
\[
\cW_{\lambda}^{\text{merge}}\colon (Y,\ve{p}\cup \{p_1,p_2\})\to (Y,\ve{p}\cup \{p_0\})
\]
for the graph cobordism which merges $p_1$ and $p_2$ into $p_0$, along the path $\lambda$. Similarly write
\[
\cW_{\lambda}^{\text{split}}\colon (Y,\ve{p}\cup \{p_0\})\to (Y,\ve{p}\cup \{p_1,p_2\})
\]
for the graph cobordism which splits $p_0$ into the pair $p_1$ and $p_2$. Write
\[
\Sp_{\lambda}:= \hat{F}_{\cW_{\lambda}^{\text{split}}}\quad \text{and} \quad \M_{\lambda}:= \hat{F}_{\cW_{\lambda}^{\text{merge}}}.
\]

\begin{figure}[ht!]
	\centering
	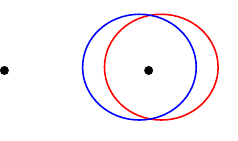
	\caption{The diagram $\cH_{p_1,(p_2)}$, considered in Lemma~\ref{lem:merge}.}\label{fig:10}
\end{figure}

\begin{lem}\label{lem:merge}
Let $\cH$ be a Heegaard diagram for $(Y,\ve{p}\cup \{p_0\})$, and let $\cH_{p_1,(p_2)}$ be the Heegaard diagram for $(Y,\ve{p}\cup \{p_1,p_2\})$ obtained by relabeling $p_0$ as $p_1$, and adding new alpha and beta curves which bound small disks containing $p_2$, as in Figure~\ref{fig:10}. Furthermore, assume that $\lambda$ is embedded in the Heegaard surface, as shown in Figure~\ref{fig:10}.  With respect to the isomorphism in equation~\eqref{eq:isomorphism-extra-basepoint}, we have
\[
\Sp_{\lambda}(\xs)=\xs\times \theta^-\quad \text{and} \quad M_{\lambda}(\xs\times \theta)=\begin{cases} \xs& \text{if } \theta=\theta^+\\
0& \text{if } \theta=\theta^-
\end{cases}.
\]
\end{lem}

\begin{figure}[ht!]
	\centering
	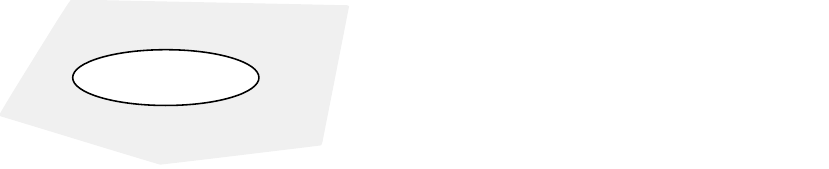
	\caption{The map for a split cobordism coincides with a contact 2-handle map. The circles represent the sutures of the manifolds $Y(\ve{p}\cup \{p_0\})$ and $Y(\ve{p}\cup \{p_1,p_2\})$.}\label{fig:3}
\end{figure}

\begin{figure}[ht!]
	\centering
	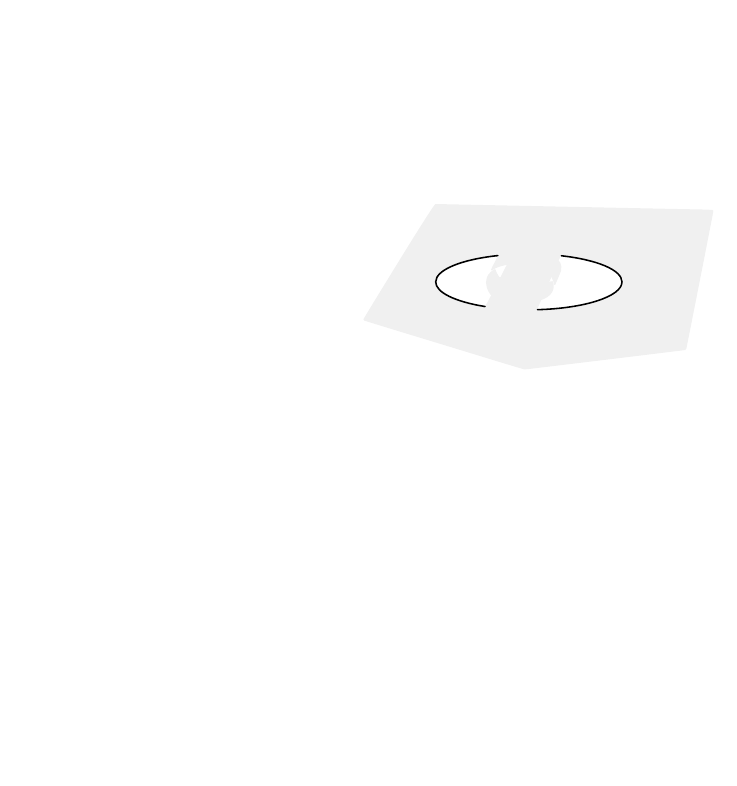
	\caption{Computing the merge map. On the left side, an index 0 holomorphic triangle is shown.}\label{fig:5}
\end{figure}

\begin{proof} We begin with the split cobordism $\cW_{\lambda}^{\text{split}}$. Write $(Z,\xi)$ for the contact portion of the boundary of the sutured manifold associated to $\cW_{\lambda}^{\text{split}}$. The contact manifold $(Z,\xi)$ is a thrice punctured, tight 3-ball. We glue $Z$ to the boundary $S^2$ of $Y(\ve{p}\cup \{p_0\})$ associated to $p_0$. The special cobordism $(\cW_{\lambda}^{\text{split}})^s$ is a product cobordism. Hence, by equation~\eqref{eq:sutured-cobordism-definition}, $\hat{F}_{\cW_{\lambda}^{\text{split}}}$ coincides with the contact gluing map $\Phi_{Z,\xi}$. The contact manifold $(Z,\xi)$ is a Morse-type contact 2-handle, so the gluing map takes the form described in Section~\ref{sec:gluing-map} (see specifically Figure~\ref{fig:9}). The description of the contact gluing map immediately gives the stated formula for $\Sp_\lambda$. See Figure~\ref{fig:3}.

We now compute the merge map. Note that the merge cobordism is obtained by turning around the split cobordism. A Morse theory argument (see \cite{JuhaszZemkeContactHandles}*{Lemma~6.7}) shows that sutured cobordism associated to $\cW^{\text{merge}}_\lambda$ has the following description:
\begin{enumerate}
\item A contact 1-handle which merges the two boundary components associated to $p_1$ and $p_2$. This turns the pair of boundary components into a single boundary component, and adds an $S^1\times S^2$ summand.
\item A 4-dimensional 2-handle which cancels the $S^1\times S^2$ summand.
\end{enumerate}
The stated formula for the merge map follows from an easy holomorphic triangle computation in the $S^1\times S^2$ summand. See Figure~\ref{fig:5}.
\end{proof}

\subsection{Spliced loop cobordisms}

We now investigate the spliced loop cobordisms, labeled \eqref{item:loop-splicing} in Figure~\ref{fig:1}.

\begin{lem}\label{lem:homology-action}
 If $\cW_\g=([0,1]\times Y, \Gamma_\g)\colon (Y,\ve{p})\to (Y,\ve{p})$ is a spliced loop cobordism, then
 \[
\hat{F}_{\cW_\g}=A_{[\g]}, 
 \]
 where $A_{[\g]}$ denotes the standard action of $H_1(Y;\Z)/\Tors$ on $\HFh(Y,\ve{p})$.
\end{lem}
\begin{proof}  Let $p\in \ve{p}$ denote the basepoint connected to the strand with the spliced loop. Write $(Z,\xi)$ for the contact portion of the boundary of the sutured manifold associated to $\cW_\g$. The contact manifold $(Z,\xi)$ has a component which consists of a twice punctured copy of $S^1\times S^2$, one of whose boundary components is glued to the boundary $S^2$ in $Y(\ve{p})$ for $p$.  The manifold $(Z,\xi)$ may be decomposed into a contact 1-handle, which splits the boundary component for $p$ into two circles (and adds no alpha or beta curves), as well as a contact 2-handle, which merges the two boundary components together, and adds an alpha and beta curve. The resulting 3-manifold is $(Y\# S^1\times S^2, \ve{p})$. Similar to Figure~\ref{fig:3}, the induced map is given by
\begin{equation}
\Phi_{Z,\xi}(\xs)= \xs\times \theta^-. \label{eq:gluing-map-spliced-loop}
\end{equation}

Let $\g_0\subset Y\# S^1\times S^2$ denote a curve which is supported in the $S^1\times S^2$ summand, and represents a generator of $H_1(S^1\times S^2)$. According to \cite{OSProperties}*{Proposition~6.4}, the map $A_{[\g_0]}$ is given by
\begin{equation}
A_{[\g_0]}(\xs\times \theta^+)=\xs\times \theta^-\quad \text{and} \quad A_{[\g_0]}(\xs\times \theta^-)=0. \label{eq:homology-action-spliced-loop}
\end{equation}

There is also a 1-handle cobordism from $(Y,\ve{\ps})$ to $(Y\# (S^1\times S^2),\ve{p})$, whose associated cobordism map is given by
\begin{equation}
F_1(\xs)=\xs\times \theta^+.\label{eq:1-handle-spliced-loop}
\end{equation}

Combining equations~\eqref{eq:gluing-map-spliced-loop}, \eqref{eq:homology-action-spliced-loop} and~\eqref{eq:1-handle-spliced-loop}, we obtain 
\[
\Phi_{Z,\xi}(\xs)=A_{[\g_0]}(F_1(\xs)).
\]

The special cobordism associated to $\cW_\g$ consists of a 2-handle, which cancels the new $S^1\times S^2$ summand. The 2-handle is attached along a framed knot $\bK$ whose underlying unframed knot is the splice $\gamma*\gamma_0$. The framing is irrelevant, since for any choice of integral framing on $\g * \g_0$, there is a canonical diffeomorphism between $(Y\# S^1\times S^2)(\bK)$ and $Y$. Hence
\begin{equation}
\hat{F}_{\cW_\g}=F_{\bK}\circ A_{[\g_0]}\circ F_{1}. \label{eq:spliced-loop}
\end{equation}
By definition, the right hand side of equation~\eqref{eq:spliced-loop} represents Ozsv\'{a}th and Szab\'{o}'s map for the identity cobordism, twisted by the induced element $[\g_0]$ of $H_1([0,1]\times Y;\Z)/\Tors$. The class in $H_1(Y;\Z)$ induced by the loop $\g_0$ coincides with $[\g]$, so the map induced by $\hat{F}_{\cW_\g}$ is exactly $A_{[\g]}$.
\end{proof}

\section{The broken path cobordism}
\label{sec:broken-path}

We now investigate the broken path cobordism, labeled \eqref{item:broken-path} in Figure~\ref{fig:1}. 
Let us write $\cB_p$ for the broken strand graph cobordism map.

 We first describe our candidate map. If $p\in \ve{p}$ and $\cH$ is a Heegaard diagram for $(Y,\ve{p})$, then there is a map
\[
\Phi_{p}\colon \CFh(\cH)\to \CFh(\cH),
\]
given by the formula
\[
\Phi_p(\xs)=\sum_{\ys\in \bT_{\a}\cap \bT_{\b}}\sum_{\substack{\phi\in \pi_2(\xs,\ys)\\ \mu(\phi)=1\\ n_p(\phi)=1\\ n_{p'}(\phi)=0, \, p'\in \ve{p}\setminus \{p\}}} \# (\cM(\phi)/\R)\cdot \ys.
\]
By counting the ends of moduli spaces of index 2 holomorphic disks which cover $p$ exactly once, we see that $\Phi_p$ is a chain map. By counting the ends of moduli spaces of index 2 holomorphic disks which cover $p$ exactly twice, we obtain
\[
\Phi_p^2=\d \circ H+H\circ \d,
\]
where $H$ is the map which counts index $1$ holomorphic disks representing classes $\phi$ with $n_{p}(\phi)=2$ and $n_{p'}(\phi)=0$ for all $p'\in \ve{p}\setminus \{p\}$.

In this section, we prove the following:

\begin{prop}\label{prop:d1-comp}
If $(Y,\ve{p})$ is a multi-pointed 3-manifold and $p\in \ve{p}$, then
 \[
\cB_p=\Phi_p,
\] 
as endomorphisms of $\HFh(Y,\ve{p})$.
\end{prop}

To prove Proposition~\ref{prop:d1-comp}, it is helpful to consider the minus version of the Heegaard Floer chain complexes. Write $\ve{p}=\{p_1,\dots, p_n\}$. We now describe an algebraic interpretation of $\Phi_{p_i}$ in terms of the chain complex $\CF^-(Y,\ve{p},\frs)$, which we recall is finitely generated and free over the ring 
\[
R_n:=\bF[U_1,\dots, U_n].
\]
  Given a Heegaard diagram $\cH$ of $(Y,\ve{p})$, the intersection points $\xs\in \bT_{\a}\cap \bT_{\b}$ with $\frs_{\ve{p}}(\xs)=\frs$ give a free basis of $\CF^-(\cH,\frs)$ over $R_n$.  The complex $\CFh(\cH,\frs)$ is obtained by setting $U_1=U_2=\cdots =U_n=0$, or equivalently by taking a tensor product with the ring $\bF$, with the trivial action of $U_i$.

 We may write the differential of $\CF^-(\cH,\frs)$ as a square matrix, using the basis of intersection points. The map $\Phi_{p_i}$ is given by taking this matrix, differentiating each entry with respect to $U_i$, and setting $U_1=\cdots=U_n=0$.

More generally, suppose $(C^-,\d^-)$ is a free, finitely generated chain complex over the ring $R_n$, with some chosen basis. Write $(\hat{C},\hat{\d})$ for the chain complex obtained by setting $U_1=\cdots=U_n=0$. We may define a map 
\[
\Phi_{U_i}\colon \hat{C}\to \hat{C},
\]
by taking the matrix for $\d^-$, and differentiating each entry with respect to $U_i$, and then setting all variables to be zero.

\begin{lem} \label{lem:commute-phi}
\begin{enumerate}
\item\label{part-1-commute-phi} Suppose $(C_1^-,\d_1^-)$ and $(C_2^-,\d_2^-)$ are free, finitely generated chain complexes over $R_n$, with fixed bases, and $F\colon C_1^-\to C_2^-$ is an $R_n$-equivariant chain map. Write $\hat{F}\colon \hat{C}_1\to \hat{C}_2$ for the induced map. Then
\[
\Phi_{U_i}\circ \hat{F}+\hat{F}\circ \Phi_{U_i}\simeq 0.
\]
\item \label{part-2-commute-phi} Suppose that $\ve{p}=\{p_1,\dots,p_n\}$ is a collection of basepoints on $Y$. The map $\Phi_p\colon \hat{\CF}(\cH)\to \hat{\CF}(\cH)$ is natural, in the sense that if $\cH$ and $\cH'$ are two diagrams for $(Y,\ve{p})$ then
\[
\Psi_{\cH\to \cH'} \circ \Phi_{p}\simeq \Phi_{p} \circ \Psi_{\cH\to \cH'},
\]
where $\Psi_{\cH\to \cH'}$ denotes the change of diagrams map from $\CFh(\cH)$ to $\CFh(\cH')$.
\end{enumerate}
\end{lem}
\begin{proof} The second statement follows from the first, since the transition maps on $\CFh$ are restrictions of the transition maps on $\CF^-$.

To prove the first claim, we take the equation
\[
0=\d^- \circ F+F\circ \d^-,
\]
and differentiate it with respect to $U_i$. Using the Leibniz rule for products of matrices, and then setting $U_1=\cdots =U_n=0$, we obtain
\[
F \circ \Phi_{U_i}+\Phi_{U_i}\circ F= \hat{\d} \circ \hat{F'}+\hat{F'}\circ \hat{\d},
\]
as maps from $\hat{C}_1$ to $\hat{C}_2$. Here, $F'$ denotes the map obtained by taking the matrix for $F$, and differentiating each entry with respect to $U_i$, and $\hat{F'}$ denotes the map resulting from setting $U_1=\cdots =U_n=0$.  
\end{proof}

\begin{proof}[Proof of Proposition~\ref{prop:d1-comp}] To disambiguate terms, let us write $\Phi_{p}^{(Y,\ve{p})}$ and $\cB_p^{(Y,\ve{p})}$ for the endomorphisms $\Phi_p$ and $\cB_p$ of $\HFh(Y,\ve{p})$.

 As a first step, we show the claim when the component of $Y$ which contains $p$ also contains another basepoint $p_0$. In this case, the $\cB_p^{(Y,\ve{p})}$ is equal to $S_p^+\circ S_p^-$, where $S_p^{\pm}$ denote the free-stabilization and destabilization maps considered in Section~\ref{sec:free-stabilization}. We may use the diagram $\cH_{(p)}$ shown in Figure~\ref{fig:2}. Using Lemma~\ref{lem:compute-free-stabilization}, we obtain that 
\begin{equation}
\cB_p^{(Y,\ve{p})}(\xs\times \theta^+)=0 \quad \text{and} \quad \cB_p^{(Y,\ve{p})}(\xs\times \theta^-)= \xs\times \theta^+. \label{eq:broken-strand-formula}
\end{equation}
On the other hand, using the diagram in Figure~\ref{fig:2}, the only holomorphic curves of index 1 going over $p$ exactly once have domain consisting of the bigon going over $p$. Using this diagram, we see that $\Phi_{p}^{(Y,\ve{p})}$ coincides with equation~\eqref{eq:broken-strand-formula}. Hence, the claim follows if there is another basepoint in the component of $Y$ which contains $p$.

We now consider the case when $p$ is the only basepoint in its component of $Y$.  In this case, we argue by adding a trivial strand to the graph, as shown in Figure~\ref{fig:6}. Adding a trivial strand does not change the isotopy class of a regular neighborhood of the graph, and hence does not change the cobordism map. We decompose the broken strand cobordism as follows:
\begin{enumerate}
\item A free-stabilization cobordism, adding a new basepoint $p_0$.
\item  A broken strand cobordism from $(Y,p,p_0)$ to $(Y,p,p_0)$ (which is broken over $p$).
\item A basepoint merging cobordism, which merges $p$ and $p_0$ along a path $\lambda$.
\end{enumerate}

\begin{figure}[ht!]
	\centering
%% Accompanies image file 'fig6.pdf' (pdf, eps, ps)
\begingroup%
  \makeatletter%
  \providecommand\color[2][]{%
    \errmessage{(Inkscape) Color is used for the text in Inkscape, but the package 'color.sty' is not loaded}%
    \renewcommand\color[2][]{}%
  }%
  \providecommand\transparent[1]{%
    \errmessage{(Inkscape) Transparency is used (non-zero) for the text in Inkscape, but the package 'transparent.sty' is not loaded}%
    \renewcommand\transparent[1]{}%
  }%
  \providecommand\rotatebox[2]{#2}%
  \newcommand*\fsize{\dimexpr\f@size pt\relax}%
  \newcommand*\lineheight[1]{\fontsize{\fsize}{#1\fsize}\selectfont}%
  \ifx\svgwidth\undefined%
    \setlength{\unitlength}{184.50007875bp}%
    \ifx\svgscale\undefined%
      \relax%
    \else%
      \setlength{\unitlength}{\unitlength * \real{\svgscale}}%
    \fi%
  \else%
    \setlength{\unitlength}{\svgwidth}%
  \fi%
  \global\let\svgwidth\undefined%
  \global\let\svgscale\undefined%
  \makeatother%
  \begin{picture}(1,0.34146329)%
    \lineheight{1}%
    \setlength\tabcolsep{0pt}%
    \put(0,0){\includegraphics[width=\unitlength,page=1]{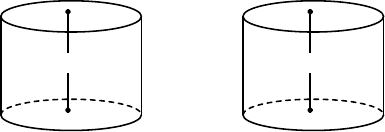}}%
    \put(0.50686514,0.1562137){\color[rgb]{0,0,0}\makebox(0,0)[t]{\lineheight{1.25}\smash{\begin{tabular}[t]{c}$=$\end{tabular}}}}%
    \put(0,0){\includegraphics[width=\unitlength,page=2]{fig6.pdf}}%
  \end{picture}%
\endgroup%

	\caption{Adding a trivial strand to the broken path cobordism map.}\label{fig:6}
\end{figure}

 Hence
\begin{equation}
\cB_p^{(Y,p)}=\M_\lambda \circ \cB_p^{(Y,p,p_0)} \circ \St_{p_0}^+.\label{eq:decompose-single-strand-broken}
\end{equation}
By the proof when there are at least 2 basepoints, equation~\eqref{eq:decompose-single-strand-broken} gives
\begin{equation}
\Phi_p^{(Y,p)}=\M_\lambda \circ \Phi_{p}^{(Y,p,p_0)} \circ \St_{p_0}^+. \label{eq:manipulate-added-stand}
\end{equation}

If we can show 
\begin{equation}
\Phi_p^{(Y,p,p_0)}\circ \St_{p_0}^+=\St_{p_0}^+\circ \Phi_p^{(Y,p)}, \label{eq:phiS^+commute}
\end{equation} 
then we can manipulate equation~\eqref{eq:manipulate-added-stand} to obtain
\begin{equation}
\begin{split}
\cB_p^{(Y,p)}&=\M_{\lambda}\circ \Phi_{p}^{(Y,p,p_0)} \circ S_{p_0}^+\\
&=\M_\lambda\circ  \St_{p_0}^+ \circ \Phi_p^{(Y,p)}\\
&=\Phi_p^{(Y,p)},
\end{split}
\label{eq:final-manipulation}
\end{equation}
since $\M_\lambda\circ \St_{p_0}^+=\id$, as the corresponding cobordism is the identity graph cobordism, with a trivial strand. Hence, it suffices to prove equation~\eqref{eq:phiS^+commute}. 

Consider the 2-variable polynomial ring $\bF[U,U_0]$, where $U$ is associated to $p$, and $U_0$ is associated to $p_0$. 
Note that $\hat{\CF}(Y,p,p_0,\frs)$ is obtained from $\CF^-(Y,p,p_0,\frs)$ by setting $U=U_0=0$. Similarly, $\CFh(Y,p,\frs)$ is also obtained from $\CF^-(Y,p)\otimes_{\bF} \bF[U_0]$ by setting $U=U_0=0$. Hence, by part~\eqref{part-1-commute-phi} of Lemma~\ref{lem:commute-phi}, to show equation~\eqref{eq:phiS^+commute}, it suffices to show that the map $\St_{p_0}^+$ can be extended to a $\bF[U,U_0]$-equivariant map from $\CF^-(Y,p,\frs)\otimes_{\bF} \bF[U_0]$ to $\CF^-(Y,p,p_0,\frs)$.

 Let $\cH$ be a diagram for $(Y,p)$, and consider a diagram $\cH_{(p_0)}$ like the one shown in  Figure~\ref{fig:2}, but with the basepoint $p_0$ encircled by the new alpha and beta circles. There is an isomorphism of modules
\[
\CF^-(\cH_{(p_0)},\frs)\iso \CF^-(\cH,\frs)\otimes_{\bF} \langle \theta^+,\theta^- \rangle \otimes_{\bF} \bF[U_0]
\]
 Ozsv\'{a}th and Szab\'{o} \cite{OSLinks}*{Equation~20} prove that there is an almost complex structure so that 
the differential on $\CF^-(\cH_{(p_0)},\frs)$ takes the form
\begin{equation}
\d_{\cH_{(p_0)}}(\xs\times \theta^+)=\d_{\cH}(\xs)\otimes \theta^+\quad \text{and}\quad  \d_{\cH_{(p_0)}}(\xs\times \theta^-)=\d_{\cH}(\xs)\otimes \theta^-   +   (U_p+U_{p_0})\cdot \xs\times \theta^+.
\label{eq:minus-differential}
\end{equation}
Equation~\eqref{eq:minus-differential} implies that the map $\xs\mapsto \xs\otimes \theta^+$, extended equivariantly over $\bF[U,U_0]$, gives an $\bF[U,U_0]$-equivariant chain map from $\CF^-(\cH,\frs)\otimes \bF[U_0]$ to $\CF^-(\cH_{(p_0)},\frs)$, which restricts to $S_{p_0}^+$ when we set $U=U_0=0$. Part~\eqref{part-1-commute-phi} of Lemma~\ref{lem:commute-phi} implies equation~\eqref{eq:phiS^+commute}, which allows us to perform the manipulation from equation~\eqref{eq:final-manipulation}, completing the proof.
\end{proof}

\section{The action of the fundamental group}

We are now ready to compute the action of the fundamental group on $\HFh(Y,\ve{p})$.

\begin{thm}\label{thm:pi-1-action} The action of $\gamma\in \pi_1(Y,p)$ on $\HFh(Y,\ve{p})$ is given by the formula
\[
\gamma_*=\id+A_{[\g]}\circ \Phi_p,
\]
where $A_{[\g]}$ denotes the action of $H_1(Y;\Z)/\Tors$.
\end{thm}

\begin{figure}[ht!]
	\centering
%% Accompanies image file 'fig7.pdf' (pdf, eps, ps)
\begingroup%
  \makeatletter%
  \providecommand\color[2][]{%
    \errmessage{(Inkscape) Color is used for the text in Inkscape, but the package 'color.sty' is not loaded}%
    \renewcommand\color[2][]{}%
  }%
  \providecommand\transparent[1]{%
    \errmessage{(Inkscape) Transparency is used (non-zero) for the text in Inkscape, but the package 'transparent.sty' is not loaded}%
    \renewcommand\transparent[1]{}%
  }%
  \providecommand\rotatebox[2]{#2}%
  \newcommand*\fsize{\dimexpr\f@size pt\relax}%
  \newcommand*\lineheight[1]{\fontsize{\fsize}{#1\fsize}\selectfont}%
  \ifx\svgwidth\undefined%
    \setlength{\unitlength}{187.13283738bp}%
    \ifx\svgscale\undefined%
      \relax%
    \else%
      \setlength{\unitlength}{\unitlength * \real{\svgscale}}%
    \fi%
  \else%
    \setlength{\unitlength}{\svgwidth}%
  \fi%
  \global\let\svgwidth\undefined%
  \global\let\svgscale\undefined%
  \makeatother%
  \begin{picture}(1,0.25717009)%
    \lineheight{1}%
    \setlength\tabcolsep{0pt}%
    \put(0,0){\includegraphics[width=\unitlength,page=1]{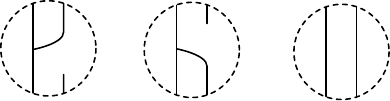}}%
    \put(0.27693503,0.12000826){\color[rgb]{0,0,0}\makebox(0,0)[lt]{\lineheight{1.25}\smash{\begin{tabular}[t]{l}$+$\end{tabular}}}}%
    \put(0.6615461,0.11899611){\color[rgb]{0,0,0}\makebox(0,0)[lt]{\lineheight{1.25}\smash{\begin{tabular}[t]{l}$=$\end{tabular}}}}%
  \end{picture}%
\endgroup%

	\caption{A local relation satisfied by the graph cobordisms.}
	\label{fig:7}
\end{figure}

As a helpful first step, we prove the relation shown in Figure~\ref{fig:7}.

\begin{lem}\label{lem:local-relation}
 The graph cobordism maps satisfy the local relation shown in Figure~\ref{fig:7}.
\end{lem}
\begin{proof}
 We view the local relation as taking place in a cylinder $[0,1]\times Y$. Let $p_1$ and $p_2$ be two basepoints of $Y$, corresponding to the two bottom points in the local relation, and let $\lambda\subset Y$  be the corresponding path connecting them. We may view the left cobordism of Figure~\ref{fig:7} as a free-destabilization, followed by a basepoint splitting cobordism. The middle cobordism is a basepoint merge, followed by a free-stabilization. The right hand side is the identity. Hence,
 it is sufficient to check
 \begin{equation}
\Sp_\lambda \circ S_{p_2}^- +S_{p_2}^+\circ M_{\lambda}=\id. \label{eq:local-relation-in-maps}
 \end{equation}
 Equation~\eqref{eq:local-relation-in-maps} is easily verified from Lemmas~\ref{lem:compute-free-stabilization} and~\ref{lem:merge}.
\end{proof}

\begin{figure}[ht!]
	\centering
%% Accompanies image file 'fig8.pdf' (pdf, eps, ps)
\begingroup%
  \makeatletter%
  \providecommand\color[2][]{%
    \errmessage{(Inkscape) Color is used for the text in Inkscape, but the package 'color.sty' is not loaded}%
    \renewcommand\color[2][]{}%
  }%
  \providecommand\transparent[1]{%
    \errmessage{(Inkscape) Transparency is used (non-zero) for the text in Inkscape, but the package 'transparent.sty' is not loaded}%
    \renewcommand\transparent[1]{}%
  }%
  \providecommand\rotatebox[2]{#2}%
  \newcommand*\fsize{\dimexpr\f@size pt\relax}%
  \newcommand*\lineheight[1]{\fontsize{\fsize}{#1\fsize}\selectfont}%
  \ifx\svgwidth\undefined%
    \setlength{\unitlength}{263.91458205bp}%
    \ifx\svgscale\undefined%
      \relax%
    \else%
      \setlength{\unitlength}{\unitlength * \real{\svgscale}}%
    \fi%
  \else%
    \setlength{\unitlength}{\svgwidth}%
  \fi%
  \global\let\svgwidth\undefined%
  \global\let\svgscale\undefined%
  \makeatother%
  \begin{picture}(1,0.28742147)%
    \lineheight{1}%
    \setlength\tabcolsep{0pt}%
    \put(0,0){\includegraphics[width=\unitlength,page=1]{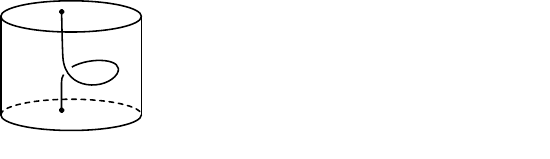}}%
    \put(0.12603571,0.19587938){\color[rgb]{0,0,0}\makebox(0,0)[lt]{\lineheight{1.25}\smash{\begin{tabular}[t]{l}$\hat{\g}$\end{tabular}}}}%
    \put(0,0){\includegraphics[width=\unitlength,page=2]{fig8.pdf}}%
    \put(0.31461085,0.16356879){\color[rgb]{0,0,0}\makebox(0,0)[t]{\lineheight{1.25}\smash{\begin{tabular}[t]{c}$=$\end{tabular}}}}%
    \put(0.68530776,0.16356879){\color[rgb]{0,0,0}\makebox(0,0)[t]{\lineheight{1.25}\smash{\begin{tabular}[t]{c}$+$\end{tabular}}}}%
    \put(0,0){\includegraphics[width=\unitlength,page=3]{fig8.pdf}}%
    \put(0.12926267,0.00313051){\color[rgb]{0,0,0}\makebox(0,0)[t]{\lineheight{1.25}\smash{\begin{tabular}[t]{c}$\g_*$\end{tabular}}}}%
    \put(0.49995935,0.00267161){\color[rgb]{0,0,0}\makebox(0,0)[t]{\lineheight{1.25}\smash{\begin{tabular}[t]{c}$\id$\end{tabular}}}}%
    \put(0.87065594,0.00313051){\color[rgb]{0,0,0}\makebox(0,0)[t]{\lineheight{1.25}\smash{\begin{tabular}[t]{c}$\Phi_p\circ A_{[\g]}$\end{tabular}}}}%
    \put(0,0){\includegraphics[width=\unitlength,page=4]{fig8.pdf}}%
  \end{picture}%
\endgroup%

	\caption{Obtaining the formula for the $\pi_1$-action by applying the local relation from Figure~\ref{fig:7} to the graph $\hat{\g}$.}
	\label{fig:8}
\end{figure}

We now prove the formula for the $\pi_1$-action:

\begin{proof}[Proof of Theorem~\ref{thm:pi-1-action}]
We focus on the case when $Y$ has a single basepoint, to simplify the notation.
The diffeomorphism map $\gamma_*$ coincides with the graph cobordism map for  $([0,1]\times Y,\hat{\g})$, where 
\[
\hat{\g}:=\{(t,\g(t)):t\in [0,1]\}\subset [0,1]\times Y.
\]
We apply the local relation from Figure~\ref{fig:7} to the graph $\hat{\g}$, as shown in Figure~\ref{fig:8}. We obtain the sum of the two graph cobordisms shown on the right side of Figure~\ref{fig:8}. We may identify the right most term with the map $\Phi_p\circ A_{[\g]}$ using Lemma~\ref{lem:homology-action} and Proposition~\ref{prop:d1-comp}. The proof is complete.
\end{proof}

\bibliographystyle{custom}
\bibliography{biblio}

\end{document}